\documentclass[12pt]{article}
\usepackage{amsmath,amssymb,amsthm,amscd,euscript}
 \DeclareMathOperator{\im}{Im}
 
\DeclareMathOperator{\SO}{SO} \DeclareMathOperator{\Or}{O}
\DeclareMathOperator{\SL}{SL} \DeclareMathOperator{\diag}{diag}

\DeclareMathOperator{\const}{const}

\def\pd#1#2{\frac{\partial{#1}}{\partial{#2}}}
\def\pd1#1{\frac{\partial}{\partial#1}}
\def\d1#1{\frac{d}{d#1}}

\def\text#1{\mbox{#1}}

\newtheorem{theore}{Theorem}\newtheorem{Lem}{Lemma}

\allowdisplaybreaks

\numberwithin{theore}{section} \numberwithin{equation}{section}
\numberwithin{proposit}{section}
\numberwithin{Def}{section}\numberwithin{Lem}{section}

\begin{document}

\author{Alexey V.~Shchepetilov\footnote{Department of Physics, Moscow State
University, 119992 Moscow, Russia, e-mail address:
quant@phys.msu.ru}}
\title{Nonintegrability of the two-body problem in constant curvature spaces II.}
\date{}\maketitle
\begin{abstract}
We consider the reduced two-body problem with a central potential
on the sphere ${\bf S}^{2}$ and the hyperbolic plane ${\bf
H}^{2}$. For two  potentials different from the Newton and the
oscillator ones we prove the nonexistence of an additional
meromorphic integral for the complexified dynamic systems. \vskip
20pt

\end{abstract}

\section{Introduction}\label{Introduction}\markright{\ref{Introduction}
Introduction}

The study of classical mechanics in constant curvature spaces was
begun in the second half of the 19th century in works of Lipschitz
\cite{Lip3}, Killing \cite{Kil2} and Neumann \cite{Neumann}. They
dealt with the one-body problem in a central potential on spheres
${\bf S}^{2}$ and ${\bf S}^{3}$. This problem revealed much
similarity with its Euclidean analogue: the existence of the
Newton-like potential (known already to Lobachevski \cite{Lob})
and three Kepler laws for it.

At the beginning of the 20th century these results were
transferred onto the hyperbolic space by Liebmann \cite{Lib02},
\cite{Lib}, who also generalized the Bertrand theorem \cite{Ber}
for constant curvature spaces \cite{Lib03}.

After the rise of special and general relativity, these results
were almost completely forgotten (see nevertheless \cite{DomZitt})
and rediscovered many times in the framework of the theory of
integrable dynamical systems (see \cite{ShchBook}, section 6.4).

The two-body problem with a central interaction in constant
curvature spaces $\mathbf{S}^{n}$ and $\mathbf{H}^{n}$
considerably differs from its Euclidean analogue. The variable
separation for the latter problem is trivial, while for the former
one no central potentials are known that admit a variable
separation.

The reduction of the two-body problem in constant curvature spaces
to a Hamiltonian system with two degrees of freedom was carried
out using different approaches in \cite{Shch98} and \cite{Sh001}
(see also \cite{ShchBook}).

The Morales-Ramis theory \cite{MorRu} contains a straightforward
method for proving the nonintegrability of Hamiltonian dynamical
systems. To do this one needs:
\begin{enumerate}
\item select a particular explicit solution of a system;
\item write a system of normal variation equations (NVE) in a neighborhood of this
solution;
\item prove that the identity component of the differential Galois
group (\cite{Kap}, \cite{vdPS}) of this NVE is not abelian.
\end{enumerate}
Note, that the direct algorithm for the third step exists only for
linear differential equations with rational coefficients.

This method (steps 1-3) was realized for many dynamical systems,
see, for example, \cite{MSS}, \cite{MacPrz} and references
therein.

The meromorphic nonintegrability of the reduced two-body problem
on the sphere ${\bf S}^{2}$ and the hyperbolic plane ${\bf H}^{2}$
with the Newton- and oscillator-like potentials was proved in
\cite{Shch06}. The particular explicit solution used there is a
motion of both bodies along a common geodesic. The extension of
this result onto the case of ${\bf S}^{3}$ and ${\bf H}^{3}$ (most
general for the two-body problem on constant curvature spaces
\cite{Shch98}) is hampered by the absence of a known explicit
solution, different from the mentioned above, in the case of
arbitrary masses.

Here, we shall prove the nonexistence of an additional meromorphic
first integral for the restricted two-body problem on
$\mathbf{S}^{2}$ and $\mathbf{H}^{2}$ with two other central
potentials, admitting the reduction of a corresponding NVE to a
second order differential equation with rational coefficients. A
more general description of a Hamiltonian reduction of systems
under consideration can be found in \cite{ShchBook}.

\section{Reduced two-body systems on $\mathbf{H}^{2},\mathbf{S}^{2}$ and NVE}\label{RTBP}
\markright{\ref{RTBP} Reduced two-body systems}

\subsection{The reduced two-body problem on the hyperbolic plane $\mathbf{H}^{2}$}

After the excluding of the diagonal (corresponding to the
collision of particles), the configuration space ${\bf
H}^{2}\times{\bf H}^{2}$ can be represented as a direct product
$$Q:=\left({\bf H}^{2}\times{\bf
H}^{2}\right)\backslash\diag\simeq I\times\Or_0(1,2),$$ where
$I:=(0,\infty)$ and the identity component $\Or_0(1,2)$ of the
isometry group for $\mathbf{H}^{2}$ acts only onto the second
factor. Therefore, the phase space
$$
M:=T^{*}Q=T^{*}I\times T^{*}\left(\Or_0(1,2)\right)
$$
is reduced to $T^{*}I\times\mathcal{O}$, where $\mathcal{O}$ is a
coadjoint orbit of $\Or_0(1,2)$.

After a time rescaling the Hamiltonian function can be represented
in the form
\begin{equation}\label{hhSimple}
h_{h}=\frac{1}{2\mu}\left(p_{\theta}^{2}+\frac{p_{2}^{2}}{\sinh^{2}\theta}\right)+
p_{\theta}p_{0}+p_{2}^{2}+p_{1}p_{2}\coth\theta+V(\theta),
\end{equation}
where $\mu:=m_{1}/(m_{1}+m_{2})\neq0$ for bodies' masses $m_{1},
m_{2}$ and the Poisson brackets for variables
$\theta,p_{\theta},p_{0},p_{1},p_{2}$ are as follows:
\begin{gather*}
\{\theta,p_{\theta}\}=1,\,\{p_{0},p_{1}\}=p_{2},\,\{p_{1},p_{2}\}=-p_{0},\\
\{p_{0},p_{2}\}=p_{1},\,\{\theta,p_{i}\}=0,\,\{p_{\theta},p_{i}\}=0,\,i=0,1,2.
\end{gather*}
Here $p_{\theta}$ is the momentum, conjugated to $\theta\in I$,
and $p_{0}^{2}+p_{1}^{2}-p_{2}^{2}=\gamma=\const\in\mathbb{R}$ on
$\mathcal{O}$.

The motion of bodies along a common geodesic with a total nonzero
momentum $\gamma$ corresponds to $p_{0}\equiv p=\const\ne 0,
p_{1}\equiv p_{2}\equiv0$ for $\gamma=p^{2}>0$. It is described by
the Hamiltonian function
\begin{equation}\label{GammaHamF}
h_{0}=\frac{1}{2\mu}p_{\theta}^{2}+
pp_{\theta}+V(\theta)=\frac{1}{2\mu}\left(p_{\theta}+\mu
p\right)^{2}+V(\theta)-\frac{\mu}2p^{2}.
\end{equation}

Let $p_{\theta}=p_{\theta}(t),\,\theta=\theta(t)$ be a solution of
the Hamiltonian system with the Hamiltonian function
(\ref{GammaHamF}). The normal variational equations in a
neighborhood of this solution are (see \cite{Shch06})
\begin{align*}\begin{split}
\frac{dp_{1}}{dt}&=-p\coth\theta
p_{1}-\left(2p+p_{\theta}+\frac{p}{\mu\sinh^{2}\theta}\right)p_{2},\\
\frac{dp_{2}}{dt}&=-p_{\theta}p_{1}+p\coth\theta p_{2}.
\end{split}
\end{align*}

\begin{theore}\label{MainThH}
The complexified Hamiltonian system with Hamiltonian function
(\ref{hhSimple}), potentials $V(\theta)=\alpha\tanh\theta$ and
$V(\theta)=\alpha\sinh^{-1}\theta,\;\alpha\in\mathbb{R}$ does not
admit an additional meromorphic first integral in the case
$m_{1}m_{2}\alpha\neq0,\,\gamma>0$.
\end{theore}
The proof of this theorem is similar to the proof of theorem
\ref{MainThS} below.

\subsection{The reduced two-body problem on the sphere $\mathbf{S}^{2}$}

To represent the configuration space ${\bf S}^{2}\times{\bf
S}^{2}$ as a direct product $I\times\SO(3)$ one should exclude not
only the diagonal, but also the set of antipodal points
$Q_{\text{op}}\cong{\bf S}^{2}$. This leads to
$$
Q_{\text{ess}}:=\left({\bf S}^{2}\times{\bf
S}^{2}\right)\backslash\left(\diag\cup Q_{\text{op}}\right)\cong
I\times\SO(3).
$$
Since
$$
\dim\left({\bf S}^{2}\times{\bf S}^{2}\right)-\dim\left(\diag\cup
Q_{\text{op}}\right)=4-2=2>1,
$$
the typical trajectory does not intersect $\diag\cup
Q_{\text{op}}$ and one can study the nonintegrability of a
Hamiltonian system on ${\bf S}^{2}\times{\bf S}^{2}$ using its
restriction onto $Q_{\text{ess}}$.

Thus, the corresponding phase space
$$
M:=T^{*}Q{\text{ess}}=T^{*}I\times T^{*}\left(\SO(3)\right)
$$
is reduced to $T^{*}I\times\mathcal{O}$, where $\mathcal{O}$ is a
coadjoint orbit of $\SO(3)$ and $I:=(0,\pi)$.

After a time rescaling the Hamiltonian function can be represented
in the form
\begin{gather}\label{hsSimple}
h_{s}=\frac{1}{2\mu}\left(p_{\theta}^{2}+\frac{p_{2}^{2}}{\sin^{2}\theta}\right)+
p_{\theta}p_{0}-p_{2}^{2}+p_{1}p_{2}\cot\theta+V(\theta),
\end{gather}
where the Poisson brackets for variables $\theta\in
I,p_{\theta},p_{0},p_{1},p_{2}$ are as follows:
\begin{gather*}
\{\theta,p_{\theta}\}=1,\,\{p_{0},p_{1}\}=-p_{2},\,\{p_{1},p_{2}\}=-p_{0},\\
\{p_{2},p_{0}\}=-p_{1},\,\{\theta,p_{i}\}=0,\,\{p_{\theta},p_{i}\}=0,\,i=0,1,2
\end{gather*}
and $p_{0}^{2}+p_{1}^{2}+p_{2}^{2}=\gamma=\const\geqslant0$ on
$\mathcal{O}$.

Again, the motion $p_{\theta}=p_{\theta}(t),\,\theta=\theta(t)$ of
bodies along a common geodesic with a total nonzero momentum
$\gamma$ corresponds to $p_{0}\equiv p=\const\ne 0, p_{1}\equiv
p_{2}\equiv0$ for $\gamma=p^{2}>0$ and is described by Hamiltonian
function (\ref{GammaHamF}).

The normal variational equations in a neighborhood of this
solution are (see \cite{Shch06})
\begin{align}\begin{split}
\frac{dp_{1}}{dt}&=-p\cot\theta
p_{1}+\left(2p+p_{\theta}-\frac{p}{\mu\sin^{2}\theta}\right)p_{2},\\
\frac{dp_{2}}{dt}&=-p_{\theta}p_{1}+p\cot\theta p_{2},
\end{split}\label{NVES}
\end{align}

Below, we shall prove
\begin{theore}\label{MainThS}
The complexified Hamiltonian system with Hamiltonian function
(\ref{hsSimple}), potentials $V(\theta)=\alpha\tan\theta$ and
$V(\theta)=-\alpha\sin^{-1}\theta,\;\alpha\in\mathbb{R}$ does not
admit an additional meromorphic first integral in the case $
m_{1}m_{2}\alpha\neq0,\,\gamma>0$.
\end{theore}

\section{Transformation of NVE's}\label{M-RT}\markright{\ref{M-RT}
Transformation of NVE's}

Here we shall transform system (\ref{NVES}) to the second order
differential equation
\begin{equation}\label{Diff2Equation2}
y''(z)=r(z)y(z).
\end{equation}

\subsection{Potential $V(\theta)=\alpha\tan\theta$}\label{ParSol1}

Consider  system (\ref{NVES}) for the potential
$V(\theta)=\alpha\tan\theta$. With respect to the new independent
variable $z:=(p_{\theta}+\mu p)/\alpha$ it can be written as
\begin{align}\begin{split}
p_{1}'(z)&=A(z)p_{1}+B(z)p_{2}\\ p_{2}'(z)&=C(z)p_{1}-A(z)p_{2}
\end{split}\label{NVESz}\end{align}
for
\begin{gather*}
A(z)=\frac{p}{f(z)(1+f^{2}(z))},\;B(z)=\frac{p}{\mu
f^{2}(z)}+\frac{(\mu-2)p-\alpha z}{1+f^{2}(z)},\\
C(z)=\frac{\alpha z-\mu p}{1+f^{2}(z)},\;
f(z)=\varepsilon-\frac{\alpha z^{2}}{2\mu},
\end{gather*}
where $\alpha\varepsilon-\mu p^{2}/2$ is a constant energy level.
Let
$$
\varkappa:=\sqrt{\frac{2\mu\varepsilon}{\alpha}},\;\eta:=\sqrt{\frac{2\mu}{\alpha}(\varepsilon-i)},\;
\lambda:=\sqrt{\frac{2\mu}{\alpha}(\varepsilon+i)},\;q:=\frac{p\mu}{\alpha};
$$
One can transform (\ref{NVESz}) into the linear differential
equation for $p_{2}(z)$ of the second order
\begin{equation}\label{p21}
p_{2}''(z)=\frac{C'}Cp'_{2}+\left(\frac{C'}CA+A^{2}+CB-A'\right)p_{2}
\end{equation}
and then into equation (\ref{Diff2Equation2}) for the function
$y(z):=p_{2}(z)/\sqrt{C}$, where
\begin{equation*}
r(z)=\frac{C'}CA+A^{2}+CB-A'-\frac12\left(\frac{C'}C\right)'+\frac14\left(\frac{C'}C\right)^{2}.
\end{equation*}
For evaluation of the function $r(z)$ one can use computer
analytical calculations, which lead to
\begin{equation}\label{r(z)}
r(z)=\sum_{j=1}^{7}\left(\frac{\alpha_{j}}{(z-z_{j})^{2}}+
\frac{\beta_{j}}{z-z_{j}}\right)=
\frac3{4z^{2}}+O(\frac1{z^{3}})\quad\text{as}\quad z\to\infty,
\end{equation}
where
$z_{1}=q,\;z_{2,3}=\pm\lambda,\;z_{4,5}=\pm\eta,\;z_{6,7}=\pm\varkappa$,
\begin{align}\label{r(z)SN}
\alpha_{1}&=\frac34,\;\alpha_{6,7}=0,\;\beta_{1}=-\frac{4q}{2q^{2}-\lambda^{2}-\eta^{2}},\;
\beta_{6}=\frac{q}{\varkappa(q-\varkappa)},\;\beta_{7}=\frac{-q}{\varkappa(q+\varkappa)},\\
\alpha_{2}&=\frac{\mu-1}{4\lambda^{2}}(q-\lambda)\left(q(\mu-1)-\lambda(\mu+1)\right),\;
\alpha_{3}=\frac{\mu-1}{4\lambda^{2}}(q+\lambda)\left(q(\mu-1)+\lambda(\mu+1)\right),\notag\\
\alpha_{4}&=\frac{\mu-1}{4\eta^{2}}(q-\eta)\left(q(\mu-1)-\eta(\mu+1)\right),\;
\alpha_{5}=\frac{\mu-1}{4\eta^{2}}(q+\eta)\left(q(\mu-1)+\eta(\mu+1)\right),\notag\\
\beta_{2}&=\frac{\mu-1}{4(\lambda^{2}-\eta^{2})\lambda^{3}}
\left(-(\mu+1)\lambda^{2}\eta^{2}+(\mu-1)\eta^{2}q^{2}-3(\mu+1)\lambda^{4}-5(\mu-1)q^{2}\lambda^{2}+8\lambda^{3}q\mu\right),\notag\\
\beta_{3}&=\frac{\mu-1}{4(\lambda^{2}-\eta^{2})\lambda^{3}}
\left((\mu+1)\lambda^{2}\eta^{2}-(\mu-1)\eta^{2}q^{2}+3(\mu+1)\lambda^{4}+5(\mu-1)q^{2}\lambda^{2}+8\lambda^{3}q\mu\right),\notag\\
\beta_{4}&=\frac{\mu-1}{4(\lambda^{2}-\eta^{2})\eta^{3}}
\left((\mu+1)\lambda^{2}\eta^{2}-(\mu-1)\lambda^{2}q^{2}+3(\mu+1)\eta^{4}+5(\mu-1)q^{2}\eta^{2}-8\eta^{3}q\mu\right),\notag\\
\beta_{5}&=\frac{\mu-1}{4(\lambda^{2}-\eta^{2})\eta^{3}}
\left(-(\mu+1)\lambda^{2}\eta^{2}+(\mu-1)\lambda^{2}q^{2}-3(\mu+1)\eta^{4}-5(\mu-1)q^{2}\eta^{2}-8\eta^{3}q\mu\right).\notag
\end{align}
Note that for $\mu=1$ the expression for $r(z)$ is considerably
simplified since this case corresponds to $m_{2}=0$ and thus to an
integrable one-body system.

\begin{Lem}\label{L1}
Suppose that
$\alpha,\varepsilon,\mu,p\in\mathbb{R},\,\mu\neq0,1,\,p\neq0$ and
\begin{equation}\label{NonEq}
(\sqrt{\varepsilon^{2}+1}-\varepsilon)(\varepsilon^{2}+1)\neq
\frac{(\mu-1)^{2}p^{2}}{4\alpha\mu};
\end{equation}
then $\alpha_{i}\notin\mathbb{R}$ for $i=2,3,4,5$.
\end{Lem}
\begin{proof}
By direct calculations one can find
$$
\im\alpha_{2}=\frac{q(\mu-1)}4\left(\frac{(1-\mu)q\alpha}{2\mu(\varepsilon^{2}+1)}\pm\frac{\sqrt{\alpha\mu}}{\sqrt{\varepsilon^{2}+1}}
\sqrt{\sqrt{\varepsilon^{2}+1}-\varepsilon}\right),
$$
therefore (\ref{NonEq}) implies $\im\alpha_{2}\ne0$.
Considerations for $\alpha_{3},\alpha_{4}$ and $\alpha_{5}$ are
similar.
\end{proof}

\subsection{Potential $V(\theta)=-\alpha\sin^{-1}\theta$}

Now, consider system (\ref{NVES}) for the potential
$V(\theta)=-\alpha\sin^{-1}\theta$. Again, w.\ r.\ t.\ the same
variable $z$, one gets
\begin{align}\begin{split}
p_{1}'(z)&=A(z)p_{1}+B(z)\sqrt{f(z)}p_{2}\\
p_{2}'(z)&=C(z)\sqrt{f(z)}p_{1}-A(z)p_{2},
\end{split}\label{NVESL}\end{align}
where
\begin{gather*}
A(z)=\frac{p}{\varphi(z)},\;B(z)=\frac{p\varphi^{2}(z)/\mu+(\mu-2)p-\alpha
z}{\varphi(z)(\varphi^{2}(z)-1)},\\ C(z)=\frac{\alpha z-\mu
p}{\varphi(z)(\varphi^{2}(z)-1)},\; \varphi(z)=\frac{\alpha
z^{2}}{2\mu}-\varepsilon,\,f(z)=\varphi^{2}(z)-1
\end{gather*}
and $\varepsilon$ is the same as above. Let
$$
\varkappa:=\sqrt{\frac{2\mu\varepsilon}{\alpha}},\;
\eta:=\sqrt{\frac{2\mu}{\alpha}(\varepsilon-1)},\;
\lambda:=\sqrt{\frac{2\mu}{\alpha}(\varepsilon+1)},\;q:=\frac{p\mu}{\alpha};
$$
then one can reduce (\ref{NVESL}) to equation
\begin{equation*}\label{p22}
p_{2}''(z)=\left(\frac{C'}C+\frac{f'}{2f}\right)p'_{2}+
\left(\left(\frac{C'}C+\frac{f'}{2f}\right)A+A^{2}+CBf-A'\right)p_{2},
\end{equation*}
which can be reduced to equation (\ref{Diff2Equation2}) by the
substitution $p_{2}(z)=y(z)\sqrt{C(z)}\left(f(z)\right)^{1/4}$.
Now it holds
\begin{gather}
z_{1}=q,\,z_{2,3}=\pm\lambda,\,z_{4,5}=\pm\eta,\,z_{6,7}=\pm\varkappa,\,\alpha_{1}=3/4,\,
\alpha_{j}=-3/16,\;j=2,3,4,5,\notag\\ \label{Lipsh}
\alpha_{6}=\frac34+(\mu-1)\left(1-\frac{q}{\varkappa}\right)
\left(1+\mu+(1-\mu)\frac{q}{\varkappa}\right),\\
\alpha_{7}=\frac34+(\mu-1)\left(1+\frac{q}{\varkappa}\right)\left(1+\mu-(1-\mu)
\frac{q}{\varkappa}\right).\notag
\end{gather}
Expressions for $\beta_{j},\,j=1,\ldots,7$ are omitted since they
will not be used below.

\begin{Lem}\label{L2}
Suppose that
$\alpha,\varepsilon,\mu,p\in\mathbb{R},\,\mu\neq0,1,\,p\neq0$ and
$\varepsilon<0$; then $\alpha_{6,7}\notin\mathbb{R}$.
\end{Lem}
\begin{proof}
One has $\varkappa\in i\mathbb{R}\backslash(0)$ and therefore
$$
i\im\alpha_{6}=-i\im\alpha_{7}=-2\mu(\mu-1)\frac{q}{\varkappa}\neq0.
$$
\end{proof}

\section{Proof of nonintegrability}\label{PNI}\markright{\ref{PNI}
Proof of nonintegrability}

In this section we shall use some facts concerning equation
(\ref{Diff2Equation2}) that were collected in the appendix of
\cite{Shch06}. For brevity, we shall cite these facts, using the
letter A.
\begin{Lem}\label{L5}
Suppose that assumptions of lemma \ref{L1} are valid. Then the
identity component $\EuScript{G}_{0}$ of the Galois group for
equation (\ref{Diff2Equation2}) with $r(z)$ given by (\ref{r(z)})
and (\ref{r(z)SN}) is not Abelian.
\end{Lem}
\begin{proof}
Here, equation (\ref{Diff2Equation2}) has eight regular singular
points $\infty,z_{i},\,i=1,\ldots,7$ and differences $\Delta_{i}$
of exponents at these points are as follows:
$\Delta_{1}=\Delta_{\infty}=2,\,\Delta_{6,7}=1,\,\Delta_{j}=\sqrt{1+4\alpha_{j}}\notin\mathbb{R}$
for $j=2,3,4,5$ (due to lemma \ref{L1}). Therefore the third case
from lemma A.1 is impossible.

Consider the first case of lemma A.1. Reasoning as in
\cite{Shch06} (the proof of lemma 5.1), one gets one of two
possibilities:
\begin{enumerate}
\item there are two linear independent solutions
$y_{k}(z),\,k=1,2$ of (\ref{Diff2Equation2}) such that
$y'_{k}/y_{k}\in\mathbb{C}(z)$ and the function $v(z)=y_{1}y_{2}$
satisfies to the equation
$$
s(z):=v'''-4rv-2r'v\equiv0;
$$
\item there is a unique (up to a constant factor) solution $y_{1}(z)$ of (\ref{Diff2Equation2}) such that
$y'_{1}/y_{1}\in\mathbb{C}(z)$ and then due to lemma A.2
$y_{1}^{m}(z)\in\mathbb{C}(z)$ for some $m\in\mathbb{Z}$ or the
identity component $\EuScript{G}_{0}$ of the differential Galois
group for (\ref{Diff2Equation2}) is nonabelian.
\end{enumerate}

For the first possibility an analysis of exponents
$\rho_{j}^{\pm}=(1\pm\Delta_{j})/2$ of equation
(\ref{Diff2Equation2}) at points $\infty,z_{j},\,j=1,\ldots,7$ and
the inclusion $v(z)\in\mathbb{C}(z)$ lead to
$$
v(z)=\frac{c}{z-q}\left(\prod_{j=2}^{5}(z-z_{j})\right),\;c=\const\neq0.
$$
Now direct calculations show that
$$
s(z)=\frac{cP_{6}(z)}{4(z-q)^{3}(z^{2}-\eta^{2})(z^{2}-\lambda^{2})(z^{2}-\varkappa^{2})^{2}},
$$
where $P_{6}(z)$ is a degree six polynomial with leading two terms
$$
-64\frac{q\mu^{2}}{\alpha^{2}}(4\mu+1)z^{6}+128\frac{q^{2}\mu^{2}}{\alpha^{2}}(4\mu-3)z^{5}
$$
and thus $s(z)\not\equiv0$. For the second possibility we
conclude, due to the complex values of exponents at
$z_{j},\,j=2,3,4,5$, that $\mathcal{G}_{0}$ is nonabelian.

Now we check if (\ref{Diff2Equation2}) has a solution of the form
$\exp\left(\int\omega(z)dz\right)$, where $\omega(z)$ is an
algebraic function over $\mathbb{C}(z)$ of degree 2. To do this we
apply the Kovacic algorithm (see the appendix in \cite{Shch06} for
notations).

One gets
$E_{z_{1}}=E_{\infty}=(-2,2,6),\,E_{z_{6}}=E_{z_{7}}=(4),\,E_{z_{j}}=(2),\,j=2,\ldots,5$.
Thus, there are no positive numbers $d(e)$ and therefore
(\ref{Diff2Equation2}) has no solutions of the form
$\exp\left(\int\omega(z)dz\right)$.

This implies that the second case of lemma A.1 can not realize.
The forth and the last case of this lemma is
$\EuScript{G}_{0}=\EuScript{G}=\SL_{2}(\mathbb{C})$. Hence, in all
possible cases the group $\EuScript{G}_{0}$ for equation
(\ref{Diff2Equation2}) with $r(z)$ given by (\ref{r(z)}) and
(\ref{r(z)SN}) is nonabelian.
\end{proof}

\begin{Lem}\label{L6}
Suppose that assumptions of lemma \ref{L2} are valid. Then the
identity component $\EuScript{G}_{0}$ of the Galois group for
equation (\ref{Diff2Equation2}) with $r(z)$ given by (\ref{r(z)})
and (\ref{Lipsh}) is not Abelian.
\end{Lem}
\begin{proof}
Now, equation (\ref{Diff2Equation2}) again has eight regular
singular points $\infty,z_{i},\,i=1,\ldots,7$ with
$\Delta_{1}=\Delta_{\infty}=2,\,\Delta_{j}=1/2,\,j=2,3,4,5,\;
\Delta_{6,7}=\sqrt{1+4\alpha_{j}}\notin\mathbb{R}$ (due to lemma
\ref{L2}). Thus, the third case from lemma A.1 is impossible.

The first case of lemma A.1 leads to one of two possibilities (see
the proof of lemma \ref{L5}). For the first possibility, the
analysis of exponents $\rho_{j}^{\pm}=(1\pm\Delta_{j})/2$ of
equation (\ref{Diff2Equation2}) at $z_{j},\,j=1,\ldots,7$ implies
$$
v(z)=\frac{P(z)}{z-q}\prod_{j=2}^{7}(z-z_{j})
$$
for a polynomial $P(z)\not\equiv0$, but since
$\rho_{\infty}\in(-1/2,3/2)$ the growth of $v(z)$ at infinity can
not be faster than $z^{3}$, which leads to the contradiction.
Thus, the first possibility can not realize and the second one
corresponds to a nonabelian $\EuScript{G}_{0}$ due to
$\Delta_{6,7}\notin\mathbb{R}$.

Check the second case of lemma A.1 using the Kovacic algorithm
(see again the appendix in \cite{Shch06}). Now
$E_{z_{1}}=E_{\infty}=(-2,2,6),\,E_{z_{6}}=E_{z_{7}}=(2),\,E_{z_{j}}=(1,2,3),\,j=2,\ldots,5$
and the maximal value for $d(e)$ equals $0$. Thus, one should
verify if the function
$$
\Theta(z)=\frac12\left(-\frac2{z-z_{1}}+\sum_{j=2}^{5}\frac{1}{z-z_{j}}+
\sum_{j=6}^{7}\frac{2}{z-z_{j}}\right)
$$
satisfies to the equation
$$
\Xi(z):=\Theta''+3\Theta\Theta'+\Theta^{3} -4r\Theta-2r'=0.
$$
But computer calculations shows that
$$
\Xi(z)=\frac{-12qz^{6}+24q^{2}z^{5}+12q\varkappa^{2}z^{4}-48q^{2}\varkappa^{2}z^{3}+\ldots}
{(z-q)^{2}(z^{2}-\lambda^{2})(z^{2}-\eta^{2})(z^{2}-\varkappa^{2})^{2}}.
$$

This implies that the second case of lemma A.1 can not realize and
as in the proof of the preceding lemma we conclude that the group
$\EuScript{G}_{0}$ for equation (\ref{Diff2Equation2}) with $r(z)$
given by (\ref{r(z)}) and (\ref{Lipsh}) is nonabelian.
\end{proof}

{\it Proof of theorem \ref{MainThS}}. Let the base field for
systems (\ref{NVESz}) and (\ref{NVESL}) be $\mathbb{C}(z)$. All
transformations, made while reducing these systems to equation
(\ref{Diff2Equation2}), are linear with coefficients from finite
extensions of $\mathbb{C}(z)$. Therefore they do not change the
identity components of the corresponding differential Galois
groups. Now theorem \ref{MainThS} follows from lemmas \ref{L5},
\ref{L6} and the Morales-Ramis theory (see section
\ref{Introduction}).\hfill$\square$

The proof of theorem \ref{MainThH} is completely similar.

\end{document}